\documentclass[a4paper,12pt,oneside]{article}
\usepackage{amsmath,amsfonts,amssymb,amsmath,latexsym,amsthm,enumerate}%,texdraw}

\newtheorem{thm}{Theorem}

\theoremstyle{definition}

\theoremstyle{plain}

\begin{document}
\title {Some identities of Frobenius-Euler polynomials arising from Frobenius-Euler basis}
\author{by \\Dae San Kim and Taekyun Kim}\date{}\maketitle
\begin{abstract}
\noindent In this paper, we give some new and interesting identities which are derived from the basis of Frobenius-Euler. Recently, Simsek et als(see $\lbrack13$\rbrack) have given some identities of $q$-analogue of Frobenius-Euler polynomials related to $q$-Bernstein polynomials. From the methods of our paper, we can also derive the results and identities of Simsek et als (cf.$\lbrack 13\rbrack$ ).
\end{abstract}
\section{Introduction}
Let $\lambda (\neq 1)\in \mathbf{C}$. As is well known, the Frobienius-Euler polynomials are defined by the generating function to be
\begin{align}\label{introeqn1}
\frac{1-\lambda}{e^{t}-\lambda}e^{xt}=e^{H(x\vert\lambda)t}=\sum_{n=0}^{\infty}H_{n}(x\vert\lambda)\frac{t^{n}}{n!}\,,
\end{align}
with the usual convention about replacing $H^{n}(x\vert\lambda$) by $H_{n}(x\vert\lambda$)\,\,(see $\lbrack 1\!\!-\!\!6\rbrack$)\,.\\
In the special case, $x\!=\!0$, $H_{n}(0\vert\lambda)=H_{n}(\lambda)$ are called the $n$-th Frobenius-Euler numbers.\\
Thus, by(\ref{introeqn1}), we get
\begin{equation}\label{introeqn2}
(H(\lambda)+1)^{n}-\lambda H_{n}(\lambda)=H_{n}(1\vert\lambda)-\lambda H_{n}(\lambda)=(1-\lambda)\delta_{0,n}\,\,,
\end{equation}
where $\delta_{0,n}$ is the Kronecker symbol.\\
From (\ref{introeqn1}), we can derive the following equation :
\begin{equation}\label{introeqn3}
H_{n}(x\vert\lambda)=(H(\lambda)+x)^{n}=\sum_{l=0}^{n}\left(\begin{array}{c} n\\l \end{array}\right)H_{n-l}(\lambda)x^{l}\,,\,\,(\text{see}\,\,\lbrack 6\!\!-\!\!13\rbrack)\,.
\end{equation}
Thus, by\,\,(\ref{introeqn3}), we easily see that the leading coefficient of $H_{n}(x\vert\lambda)$ is $H_{0}(\lambda)=1$\,. So, $H_{n}(x\vert\lambda)$ is a monic polynomials of degree $n$ with coefficients in $\mathbf{Q}(\lambda)$.\\ From (\ref{introeqn1}), we have
\begin{equation}\label{introeqn4}
\sum_{n=0}^{\infty}(H_{n}(x+1\vert\lambda)-\lambda H_{n}(x\vert\lambda))\frac{t^{n}}{n!}=\frac{(1-\lambda)e^{(x+1)t}}{e^{t}-\lambda}-\lambda\frac{1-\lambda}{e^{t}-\lambda}e^{xt}\,.
\end{equation}
Thus, by (\ref{introeqn4}), we get
\begin{equation}\label{introeqn5}
H_{n}(x+1\vert\lambda)-\lambda H_{n}(x\vert\lambda) = (1-\lambda)x^{n}\,,\,\,\, \text{for}\,\, n\in\mathbf{Z}_{+}\,.
\end{equation}
It is easy to show that
\begin{equation}\label{introeqn6}
\frac{d}{dx}H_{n}(x\vert\lambda) = \frac{d}{dx}(H(\lambda)+x)^{n} = n H_{n-1}(x\vert\lambda)\,,\,\,(n\in\mathbf{N})\,.
\end{equation}
From (\ref{introeqn6}), we have
\begin{equation}\label{introeqn7}
 \int_{0}^{1}H_{n}(x\vert\lambda)dx = \frac{1}{n+1}(H_{n+1}(1\vert\lambda)-H_{n+1}(\lambda))=\frac{\lambda-1}{n+1}H_{n+1}(\lambda)\,.
\end{equation}
Let $\mathbb{P}_{n}(\lambda)=\lbrace p(x)\in\mathbf{Q}(\lambda)\lbrack x\rbrack\vert$ deg $p(x)\leq n\rbrace$ be a vector space over $\mathbf{Q}(\lambda)$. Then we note that $\lbrace H_{0}(x\vert\lambda),H_{1}(x\vert\lambda),\cdots,H_{n}(x\vert\lambda)\rbrace$ is a good basis for $\mathbb{P}_{n}(\lambda)$.\\
In this paper, we develop some new methods to obtain some new identities and properties of Frobenius-Euler polynomials which are derived from the basis of the Frobenius-Euler polynomials. Those methods are useful in studying the identities of Frobenius-Euler polynomials.

\section{Some identities of Frobenius-Euler polynomials}
Let us take $p(x)\in\mathbb{P}_{n}(\lambda)$. Then $p(x)$ can be expressed as a $\mathbf{Q}(\lambda)$-linear combination of $H_{0}(x\vert\lambda),\cdots,H_{n}(x\vert\lambda)$ as follows :
\begin{equation}\label{introeqn8}
p(x) = b_{0}H_{0}(x\vert\lambda)+b_{1}H_{1}(x\vert\lambda)+\cdots+b_{n}H_{n}(x\vert\lambda)=\sum_{k=0}^{n}b_{k}H_{k}(x\vert\lambda)\,.
\end{equation}
Let us define the operator $\triangle_{\lambda}$ by
\begin{equation}\label{introeqn9}
g(x)=\triangle_{\lambda}p(x)=p(x+1)-\lambda p(x)\,.
\end{equation}
From (\ref{introeqn9}), we can derive the following equation (\ref{introeqn10}) :
\begin{equation}\label{introeqn10}
g(x)=\triangle_{\lambda}p(x)=\sum_{k=0}^{n}b_{k}( H_{k}(x+1\vert\lambda)-\lambda H_{k}(x\vert\lambda))=(1-\lambda)\sum_{k=0}^{n}b_{k}x^{k}\,.
\end{equation}
For $r\in\mathbf{Z}_{+}$, let us take the $r$-th derivative of $g(x)$ in (\ref{introeqn10}) as follows :
\begin{equation}\label{introeqn11}
g^{(r)}(x)=(1-\lambda)\sum_{k=r}^{n}k(k-1)\cdots(k-r+1)b_{k}x^{k-r}, \text{where}\,\,\, g^{(r)}(x)=\frac{d^{r}g(x)}{dx^{r}}\,.
\end{equation}
Thus, by (\ref{introeqn11}), we get
\begin{equation}\label{introeqn12}
g^{r}(0)=\frac{d^{r}g(x)}{dx^{r}}\vert_{x=0}=(1-\lambda)r!b_{r}\,.
\end{equation}
From (\ref{introeqn12}), we have
\begin{equation}\label{introeqn13}
b_{r}=\frac{g^{(r)}(0)}{(1-\lambda)r!}=\frac{1}{(1-\lambda)r!}(p^{(r)}(1)-\lambda p^{(r)}(0))\,,
\end{equation}
where $r\in\mathbf{Z}_{+}$\,, and $p^{(r)}(0)=\frac{d^{r}p(x)}{dx^{r}}\vert_{x=0}$\,.
Therefore, by (\ref{introeqn13}), we obtain the following theorem.
\begin{thm}\label{theorem1}
For $\lambda(\neq1)\in\mathbf{C}$, $n\in\mathbf{Z}_{+}$,\\
let $p(x)\in\mathbb{P}_{n}(\lambda)$ \text{with}
$p(x)=\sum_{k=0}^{n}b_{k}H_{k}(x\vert\lambda)$\,.\,\, Then we have
\begin{equation*}
b_{k}=\frac{1}{(1-\lambda)k!}g^{(k)}(0)=\frac{1}{(1-\lambda)k!}(p^{(k)}(1)-\lambda p^{(k)}(0))\,.
\end{equation*}
\end{thm}
\noindent Let us take $p(x)=H_{n}(x\vert\lambda^{-1})$. Then, by Theorem \ref{theorem1}, we get
\begin{equation}\label{introeqn14}
H_{n}(x\vert\lambda^{-1})=\sum_{k=0}^{n}b_{k}H_{k}(x\vert\lambda),
\end{equation}
where
\begin{align}\label{introeqn15}
b_{k}&= \frac{1}{(1-\lambda)k!}\frac{n!}{(n-k)!}\lbrace H_{n-k}(1\vert\lambda^{-1})-\lambda H_{n-k}(\lambda^{-1})\rbrace\\
     &=\frac{1}{1-\lambda}\left(\begin{array}{c}n\\k\end{array}\right)\lbrace H_{n-k}(1\vert\lambda^{-1})-\lambda H_{n-k}(\lambda^{-1})\rbrace\nonumber\\
     &=\frac{1}{1-\lambda}\left(\begin{array}{c}n\\k\end{array}\right)\lbrace (1-\lambda^{-1})0^{n-k}+\frac{1}{\lambda}H_{n-k}(\lambda^{-1})-\lambda H_{n-k}(\lambda^{-1})\rbrace\nonumber\,.
\end{align}
By (\ref{introeqn14}) and (\ref{introeqn15}), we get
\begin{align}\label{introeqn16}
H_{n}&(x\vert\lambda^{-1})\\
&=-\frac{1}{\lambda}H_{n}(x\vert\lambda)+\sum_{k=0}^{n}\lbrace\frac{\left(\begin{array}{c}n\\k\end{array}\right)}{\lambda(1-\lambda)}H_{n-k}(\lambda^{-1})-\frac{\lambda\left(\begin{array}{c}n\\k\end{array}\right)}{1-\lambda}H_{n-k}(\lambda^{-1})H_{k}\rbrace(x\vert\lambda)\nonumber\\
&=-\frac{1}{\lambda}H_{n}(x\vert\lambda)+\sum_{k=0}^{n}\left(\begin{array}{c}n\\k\end{array}\right)\frac{1+\lambda}{\lambda}H_{n-k}(\lambda^{-1})H_{k}(x\vert\lambda)\,.\nonumber
\end{align}
Therefore, by (\ref{introeqn16}), we obtain the following theorem.
\begin{thm}\label{theorem2}
For $n\in\mathbf{Z}_{+}$\,, we have
\begin{equation*}
\lambda H_{n}(x\vert\lambda^{-1})+H_{n}(x\vert\lambda)=(1+\lambda)\sum_{k=0}^{n}\left(\begin{array}{c}n\\k\end{array}\right)H_{n-k}(\lambda^{-1})H_{k}(x\vert\lambda)\,.
\end{equation*}
\end{thm}
\noindent Let
\begin{equation}\label{introeqn17}
p(x)=\sum_{k=0}^{n}H_{k}(x\vert\lambda)H_{n-k}(x\vert\lambda)\in\mathbb{P}_{n}(\lambda)\,.
\end{equation}
From Theorem \ref{theorem2},\,\, we note that $p(x)$ can be generated by $\lbrace H_{0}(x\vert\lambda),H_{1}(x\vert\lambda),\\\cdots,H_{n}(x\vert\lambda)\rbrace$ as follows:
\begin{equation}\label{introeqn18}
p(x)=\sum_{k=0}^{n}H_{k}(x\vert\lambda)H_{n-k}(x\vert\lambda)=\sum_{k-0}^{n}b_{k}H_{k}(x\vert\lambda)\,.
\end{equation}
By (\ref{introeqn17}), we get
\begin{equation}\label{introeqn19}
p^{(k)}(x)=\frac{(n+1)!}{(n-k+1)!}\sum_{l=k}^{n}H_{l-k}(x\vert\lambda)H_{n-k}(x\vert\lambda)\,,
\end{equation}
and
\begin{align}\label{introeqn20}
b_{k}&=\frac{1}{(1-\lambda) k!}\lbrace p^{(k)}(1)-\lambda p^{(k)}(0)\rbrace\\
     &=\frac{(n+1)!}{(1-\lambda)k!(n-k+1)!}\sum_{l=k}^{n}\lbrace H_{l-k}(1\vert\lambda)H_{n-l}(1\vert\lambda)-\lambda H_{l-k}(\lambda)H_{n-l}(\lambda)\rbrace\nonumber\\
     &=\frac{n+1}{(1-\lambda)(n-k+1)}\left(\begin{array}{c}n\\k\end{array}\right)\sum_{l=k}^{n}\lbrace(\lambda H_{l-k}(\lambda)+(1-\lambda)\delta_{0,l-k})(\lambda H_{n-l}+\nonumber\\
     &\quad(1-\lambda)\delta_{0,n-l})-\lambda H_{l-k}(\lambda)H_{n-l}(\lambda)\rbrace\nonumber\\
     &=\frac{n+1}{(1-\lambda)(n-k+1)}\left(\begin{array}{c}n\\k\end{array}\right)\sum_{l=k}^{n}\lbrace \lambda(1-\lambda)\delta_{0,l-k}H_{n-l}(\lambda)+\lambda(1-\lambda)\nonumber\\
     &\quad \times H_{l-k}(\lambda)\delta_{0,n-l}+(1-\lambda)^{2}\delta_{0,l-k}\delta_{0,n-l}+\lambda(\lambda-1)H_{l-k}(\lambda)H_{n-l}(\lambda)\rbrace\nonumber\\
     &=\frac{n+1}{(1-\lambda)(n-k+1)}\left(\begin{array}{c}n\\k\end{array}\right)\sum_{l=k}^{n}\lbrace\lambda(\lambda-1)H_{l-k}(\lambda)H_{n-l}(\lambda)+\lambda(1-\lambda)\nonumber\\
     &\quad \times H_{n-k}(\lambda)+\lambda(1-\lambda)H_{n-k}(\lambda)+(1-\lambda)^{2}\delta_{n,k}\rbrace\nonumber\\
     &=\frac{n+1}{n-k+1}\left(\begin{array}{c}n\\k\end{array}\right)\sum_{l=k}^{n}\lbrace-\lambda H_{l-k}(\lambda)H_{n-l}(\lambda)+2\lambda H_{n-k}(\lambda)+(1-\lambda)\delta_{n,k}\rbrace\,.\nonumber
\end{align}
From (\ref{introeqn18}) and (\ref{introeqn20}), we have
\begin{align}\label{introeqn21}
&\sum_{k=0}^{n}H_{k}(x\vert\lambda)H_{n-k}(x\vert\lambda)=(n+1)\sum_{k=0}^{n-1}\frac{\left(\begin{array}{c}n\\k\end{array}\right)}{n-k+1}\sum_{l=k}^{n}\lbrace(-\lambda)H_{l-k}(\lambda)H_{n-l}(\lambda)\\
&+2\lambda H_{n-k}(\lambda)\rbrace H_{k}(x\vert\lambda)+(n+1)H_{n}(x\vert\lambda)\,.\nonumber
\end{align}
Therefore, by (\ref{introeqn21}), we obtain the following theorem
\begin{thm}\label{theorem3}
For $n\in\mathbf{Z}_{+}$\,, we have
\begin{align*}
&\frac{1}{n+1}\sum_{k=0}^{n}H_{k}(x\vert\lambda)H_{n-k}(x\vert\lambda)\\
&=\sum_{k=0}^{n-1}\frac{\left(\begin{array}{c}n\\k\end{array}\right)}{n-k+1}\sum_{l=k}^{n}\lbrace(-\lambda)H_{l-k}(\lambda)H_{n-l}(\lambda)+2\lambda H_{n-k}(\lambda)\rbrace H_{k}(x\vert\lambda)+H_{n}(x\vert\lambda)\,.
\end{align*}
\end{thm}
\noindent Let us consider
\begin{equation}\label{introeqn22}
p(x)=\sum_{k=0}^{n}\frac{1}{k!(n-k)!}H_{k}(x\vert\lambda)H_{n-k}(x\vert\lambda)\in\mathbb{P}_{n}(\lambda)\,.
\end{equation}
By Theorem \ref{theorem1}, $p(x)$ can be expressed by
\begin{equation}\label{introeqn23}
p(x)=\sum_{k=0}^{n}b_{k}H_{k}(x\vert\lambda)\,.
\end{equation}
From (\ref{introeqn22}), we have
\begin{equation}\label{introeqn24}
p^{(r)}(x)=2^{r}\sum_{k=r}^{n}\frac{H_{k-r}(x\vert\lambda)H_{n-k}(x\vert\lambda)}{(k-r)!(n-k)!}\,, (r\in\mathbf{Z}_{+})\,.
\end{equation}
By Theorem \ref{theorem1}, we get
\begin{align}\label{introeqn25}
b_{k}&=\frac{1}{2k!}\lbrace p^{(k)}(1)-p^{(k)}(0)\rbrace\\
     &=\frac{2^{k-1}}{k!}\sum_{l=k}^{n}\frac{1}{(l-k)!(n-l)!}\lbrace H_{l-k}(1\vert\lambda)H_{n-l}(1\vert\lambda)-\lambda H_{l-k}(\lambda)H_{n-l}(\lambda)\rbrace\nonumber\\
     &=\frac{2^{k-1}}{k!}\sum_{l=k}^{n}\frac{1}{(l-k)!(n-l)!}\lbrace(\lambda H_{l-k}(\lambda)+(1-\lambda)\delta_{0,l-k})(\lambda H_{n-l}(\lambda)\nonumber\\
     &\quad+(1-\lambda)\delta_{0,n-l})-\lambda H_{l-k}(\lambda)H_{n-l}(\lambda)\rbrace\nonumber\\
     &=\frac{2^{k-1}}{k!}\lbrace\sum_{l=k}^{n}\frac{\lambda(\lambda-1)H_{l-k}(\lambda)H_{n-l}(\lambda)}{(l-k)!(n-l)!}+\frac{2\lambda(1-\lambda)H_{n-k}(\lambda)}{(n-k)!}+(1-\lambda)^{2}\delta_{n,k}\rbrace\nonumber\\
     &=\left\{\begin{array}{ll} \frac{2^{k-1}}{k!}\sum_{l=k}^{n}\lbrace\frac{\lambda(\lambda-1)H_{l-k}(\lambda)H_{n-l}(\lambda)}{(l-k)!(n-l)!}+\frac{2\lambda(1-\lambda)H_{n-k}(\lambda)}{(n-k)!}\rbrace\,, & \textrm{if $k\neq n$}\\
     \frac{2^{n-1}(1-\lambda)}{n!}\,,&\textrm{if $k=n$}\end{array}\right.\nonumber
\end{align}
Therefore, by (\ref{introeqn25}), we obtain the following theorem\,.
\begin{thm}\label{theorem4}
For $n\in\mathbf{Z}_{+}$, we have
\begin{align*}
\sum_{k=0}^{n}&\frac{1}{k!(n-k)!}H_{k}(x\vert\lambda)H_{n-k}(x\vert\lambda)\\
&=\sum_{k=0}^{n-1}\frac{2^{k-1}}{k!}\sum_{l=k}^{n}\lbrace\frac{\lambda(\lambda-1)H_{l-k}(\lambda)H_{n-l}(\lambda)}{(l-k)!(n-l)!}+\frac{2\lambda(1-\lambda)H_{n-k}(\lambda)}{(n-k)!}\rbrace H_{k}(x\vert\lambda)\nonumber\\
&\quad+\frac{2^{n-1}(1-\lambda)}{n!}H_{n}(x\vert\lambda)\,.
\end{align*}
\end{thm}
\section{Higher-order Frobenius-Euler polynomials}
For $n\in\mathbf{Z}_{+}$\,, the Frobenius-Euler polynomials of order $r$ are defined by the generating function to be
\begin{align}\label{introeqn26}
(\frac{1-\lambda}{e^{t}-\lambda})^{r}e^{xt}&=e^{H^{(r)}(x\vert\lambda)t}\\
&=\sum_{n=0}^{\infty}H_{n}^{(r)}(x\vert\lambda)\frac{t^{n}}{n!}\,,\nonumber
\end{align}
with the usual convention about replacing $(H^{(r)}(x\vert\lambda))^{n}$ by $ H_{n}^{(r)}(x\vert\lambda)$, (see $\lbrack1\!\!-\!\!10\rbrack$)\,.
In the special case, $x\!\!=\!\!0$, $H_{n}^{(r)}(0\vert\lambda)=H_{n}^{(r)}(\lambda)$ are called the $n$-th Frobenius-Euler numbers of order $r$\,,\,\,(see $\lbrack8\!\!-\!\!9\rbrack$)\,.\\
From (\ref{introeqn26}), we have
\begin{equation}\label{introeqn27}
H_{n}^{(r)}(x\vert\lambda)=H^{(r)}(\lambda)+x)^{n}=\sum_{l=0}^{n}\left(\begin{array}{c}n\\l\end{array}\right)H_{n-l}^{(r)}(\lambda)x^{l}\,,
\end{equation}
with the usual convention about replacing $(H^{(r)}(\lambda))^{n}$ by $H_{n}^{(r)}(\lambda)$\,.\\
By (\ref{introeqn26}), we get
\begin{equation}\label{introeqn28}
H_{n}^{(r)}(\lambda) = \sum_{n_{1}+\cdots+n_{r}=n}\left(\begin{array}{c}n\\n_{1},n_{2},\cdots,n_{r}\end{array}\right)H_{n_{1}}(\lambda)\cdots H_{n_{r}}(\lambda)\,,
\end{equation}
where $\left(\begin{array}{c}n\\n_{1},n_{2},\cdots,n_{r}\end{array}\right) = \frac{n!}{n_{1}!n_{2}!\cdots n_{r}!}$\,.
From (\ref{introeqn27}) and (\ref{introeqn28}), we note that the leading coefficient of $H_{n}^{(r)}(x\vert\lambda)$ is given by
\begin{align}\label{introeqn29}
H_{0}^{(r)}(\lambda)&=\sum_{n_{1}+\cdots+n_{r}=0}\left(\begin{array}{c}n\\n_{1},n_{2},\cdots,n_{r}\end{array}\right)H_{n_{1}}(\lambda)\cdots H_{n_{r}}(\lambda)\nonumber\\
                    &=H_{0}(\lambda)\cdots H_{0}(\lambda)=1\,.
\end{align}
Thus, by (\ref{introeqn29}), we see that $H_{n}^{(r)}$ is a monic polynomial of degree $n$ with coefficients in $\mathbf{Q}(\lambda)$\,.
From (\ref{introeqn26}), we have
\begin{equation}\label{introeqn30}
H_{n}^{(0)}(x\vert\lambda)=x^{n},\,\, \text{for}\,\, n\in\mathbf{Z}_{+}\,,
\end{equation}
and
\begin{equation}\label{introeqn31}
\frac{\partial}{\partial x}H_{n}^{(r)}(x\vert\lambda)=\frac{\partial}{\partial x}(H^{(r)}(\lambda)+x)^{n}= nH_{n-1}^{(r)}(x\vert\lambda),\,\, (r\geq 0)\,.
\end{equation}
It is not difficult to show that
\begin{equation}\label{introeqn32}
H_{n}^{(r)}(x+1\vert\lambda)-\lambda H_{n}^{(r)}(x\vert\lambda)=(1-\lambda)H_{n}^{(r-1)}(x\vert\lambda)\,.
\end{equation}
Now, we note that $\lbrace H_{0}^{(r)}(x\vert\lambda), H_{1}^{(r)}(x\vert\lambda),\cdots, H_{n}^{(r)}(x\vert\lambda)\rbrace$ is also a good basis for $\mathbb{P}_{n}(\lambda)$\,.\\
Let us define the operator $D$ as $Df(x)=\frac{df(x)}{dx}$ and let $p(x)\in\mathbb{P}_{n}(\lambda)$\,. Then  $p(x)$ can be written as
\begin{equation}\label{introeqn33}
p(x)=\sum_{k=0}^{n}C_{k}H_{k}^{(r)}(x\vert\lambda)\,.
\end{equation}
From (\ref{introeqn9}) and (\ref{introeqn32}), we have
\begin{equation}\label{introeqn34}
\triangle_{\lambda}H_{n}^{(r)}(x\vert\lambda)=H_{n}^{(r)}(x+1\vert\lambda)-\lambda H_{n}^{(r)}(x\vert\lambda)=(1-\lambda)H_{n}^{(r-1)}(x\vert\lambda)\,.
\end{equation}
Thus, by (\ref{introeqn33}) and (\ref{introeqn34}), we get
\begin{equation}\label{introeqn35}
\triangle_{\lambda}^{r}p(x)=(1-\lambda)^{r}\sum_{k=0}^{n}C_{k}H_{k}^{(0)}(x\vert\lambda)=(1-\lambda)^{r}\sum_{k=0}^{n}C_{k}x^{k}\,.
\end{equation}
Let us take the $k$-th derivative of $\triangle_{\lambda}^{r}p(x)$ in (\ref{introeqn35})\,.\\
Then we have
\begin{equation}\label{introeqn36}
D^{k}(\triangle_{\lambda}^{r}p(x))=(1-\lambda)^{r}\sum_{l=k}^{n}\frac{l!}{(l-k)!}C_{l}x^{l-k}\,.
\end{equation}
Thus, from (\ref{introeqn36}), we have
\begin{equation}\label{introeqn37}
D^{k}(\triangle_{\lambda}^{r}p(0))=(1-\lambda)^{r}\sum_{l=k}^{n}\frac{l!C_{l}}{(l-k)!}0^{l-k}=(1-\lambda)^{r}k!C_{k}\,.
\end{equation}
Thus, by(\ref{introeqn37}), we get
\begin{align}\label{introeqn38}
C_{k}&=\frac{D^{k}(\triangle_{\lambda}^{r}p(0))}{(1-\lambda)^{r}k!}\\
     &=\frac{\triangle_{\lambda}^{r}(D^{k}p(0))}{(1-\lambda)^{r}k!}=\frac{1}{(1-\lambda)^{r}k!}\sum_{j=0}^{r}\left(\begin{array}{c}r\\j\end{array}\right)(-\lambda)^{(r-j)}D^{k}p(j)\nonumber\,.
\end{align}
Therefore, by (\ref{introeqn33}) and (\ref{introeqn38}), we obtain the following theorem.
\begin{thm}\label{theorem5}
For $r\in\mathbf{Z}_{+}$\,,\,\, let $p(x)\in\mathbb{P}_{n}(\lambda)$ with
\begin{equation*}
p(x)=\frac{1}{(1-\lambda)^{r}}\sum_{k=0}^{n}C_{k}H_{k}^{(r)}(x\vert\lambda),\,\,(C_{k}\in\mathbf{Q}(\lambda))\,.
\end{equation*}
Then we have
\begin{equation*}
C_{k}=\frac{1}{(1-\lambda)^{r}k!}\sum_{j=0}^{r}\left(\begin{array}{c}r\\j\end{array}\right)(-\lambda)^{r-j}D^{k}p(j)\,.
\end{equation*}
That is,
\begin{equation*}
p(x)=\frac{1}{(1-\lambda)^{r}}\sum_{k=0}^{n}(\sum_{j=0}^{r}\frac{1}{k!}\left(\begin{array}{c}r\\j\end{array}\right)(-\lambda)^{r-j}D^{k}p(j))H_{k}^{(r)}(x\vert\lambda)\,.
\end{equation*}
\end{thm}
\noindent Let us take $p(x)=H_{n}(x\vert\lambda)\in\mathbf{P}_{n}(\lambda)$\,. Then, by Theorem \ref{theorem5}, $p(x)=H_{n}(x\vert\lambda)$ can be generated by $\lbrace H_{0}^{(r)}(x\vert\lambda), H_{1}^{(r)}(\lambda),\cdots,H_{n}^{(r)}(x\vert\lambda)\rbrace$ as follows :
\begin{equation}\label{introeqn39}
H_{n}(x\vert\lambda)=\sum_{k=0}^{n}C_{k}H_{k}^{(r)}(x\vert\lambda)\,,
\end{equation}
where
\begin{equation}\label{introeqn40}
C_{k}=\frac{1}{(1-\lambda)^{r}}\frac{1}{k!}\sum_{j=0}^{r}\left(\begin{array}{c}r\\j\end{array}\right)(-\lambda)^{r-j}D^{k}p(j)\,,
\end{equation}
and
\begin{equation}\label{introeqn41}
p^{(k)}(x)=D^{k}p(x)=n(n-1)\cdots(n-k+1)H_{n-k}(x\vert\lambda)=\frac{n!}{(n-k)!}H_{n-k}(x\vert\lambda)\,.
\end{equation}
By (\ref{introeqn40}) and (\ref{introeqn41}), we get
\begin{equation}\label{introeqn42}
C_{k}=\frac{1}{(1-\lambda)^{r}}\left(\begin{array}{c}n\\k\end{array}\right)\sum_{j=0}^{r}\left(\begin{array}{c}r\\j\end{array}\right)(-\lambda)^{r-j}H_{n-k}(j\vert\lambda)\,.
\end{equation}
Therefore, by\,(\ref{introeqn39}) and (\ref{introeqn42}), we obtain the following theorem.
\begin{thm}\label{theorem6}
For $n\in\mathbf{Z}_{+}$, we have
\begin{equation*}
H_{n}(x\vert\lambda)=\frac{1}{(1-\lambda)^{r}}\sum_{k=0}^{n}\left(\begin{array}{c}n\\k\end{array}\right)(\sum_{j=0}^{r}\left(\begin{array}{c}r\\j\end{array}\right)(-\lambda)^{r-j}H_{n-k}(j\vert\lambda))H_{k}^{(r)}(x\vert\lambda)\,.
\end{equation*}
\end{thm}
\noindent Let us assume that $p(x)=H_{n}^{(r)}(x\vert\lambda)$\,.\\
Then we have
\begin{align}\label{introeqn43}
p^{k}(x)&=n(n-1)\cdots(n-k+1)H_{n-k}^{(r)}(x\vert\lambda)\\
        &=\frac{n!}{(n-k)!}H_{n-k}^{(r)}(x\vert\lambda)\,.\nonumber
\end{align}
From Theorem \ref{theorem1}, we note that $p(x)=H_{n}^{(r)}(x\vert\lambda)$ can be expressed as a linear combination of $ H_{0}(x\vert\lambda),H_{1}(x\vert\lambda),\cdots,H_{n}(x\vert\lambda)$ :
\begin{equation}\label{introeqn44}
H_{n}^{(r)}(x\vert\lambda)=\sum_{k=0}^{n}b_{k}H_{k}(x\vert\lambda)\,,
\end{equation}
where
\begin{align}\label{introeqn45}
b_{k}&=\frac{1}{(1-\lambda)k!}\lbrace p^{k}(1)-\lambda p^{(k)}(0)\rbrace\\
     &=\frac{n!}{(1-\lambda)k!(n-k)!}\lbrace H_{n-k}^{(r)}(1\vert\lambda)-\lambda H_{n-k}^{(r)}(\lambda)\rbrace\,.\nonumber
\end{align}
By (\ref{introeqn34}) and (\ref{introeqn45}), we get
\begin{equation}\label{introeqn46}
b_{k}=\frac{1}{(1-\lambda)^{2}}\left(\begin{array}{c}n\\k\end{array}\right)H_{n-k}^{(r-1)}(\lambda)\,.
\end{equation}
Therefore, by (\ref{introeqn44}) and (\ref{introeqn46}), we obtain the following theorem.
\begin{thm}\label{theorem7}
For $n\in\mathbf{Z}_{+}$, we have
\begin{equation*}
H_{n}^{(r)}(x\vert\lambda)=\frac{1}{(1-\lambda)^{2}}\sum_{k=0}^{n}\left(\begin{array}{c}n\\k\end{array}\right)H_{n-k}^{(r-1)}(\lambda)H_{k}(x\vert\lambda)\,.
\end{equation*}
\end{thm}

\author{Department of Mathematics, Sogang University, Seoul 121-742, Republic of Korea
\\e-mail: dskim@sogang.ac.kr}\\
\\
\author{Department of Mathematics, Kwangwoon University, Seoul 139-701, Republic of Korea
\\e-mail: tkkim@kw.ac.kr}
\end{document}